\documentclass[11pt,a4paper]{article}
\usepackage{epsfig}
\usepackage{graphicx}
\usepackage[latin1]{inputenc}
\usepackage{amsmath, amsthm,  amsfonts, amscd}
\usepackage{amssymb} 
\usepackage{url}
\newtheorem{lemma}{Lemma}
\newtheorem{theorem}[lemma]{Theorem}

\newtheorem{conjecture}[lemma]{Conjecture}
\numberwithin{equation}{section} \numberwithin{lemma}{section}

\newcommand{\imply}{\Longrightarrow}

\newcommand{\sym}{\mbox{sym}}
\newcommand{\modd}{\mbox{~~mod~}}
\newcommand{\RRe}{\rm Re}
\newcommand{\IIm}{\rm Im}
\newcommand{\kommentar}[1]{}
 % ...
 \kommentar{%
   Kommentartext, wird als Argument gelesen und ignoriert
  }%

\def\thebibliography#1{\section{\centerline{\sc References}}
  \global\def\@listi{\leftmargin\leftmargini
               \labelwidth\leftmargini \advance\labelwidth-\labelsep
               \topsep 1pt plus 2pt minus 1pt
               \parsep 0.25ex plus 1pt \itemsep 0.25ex plus 1pt}
  \list {[\arabic{enumi}]}{\settowidth\labelwidth{[#1]}\leftmargin\labelwidth
    \advance\leftmargin\labelsep\usecounter{enumi}}
    \def\newblock{\hskip .11em plus .33em minus -.07em}
    \sloppy
    \sfcode`\.=1000\relax}

\begin{document}

\title{\bf{LOWER ORDER TERMS FOR THE ONE-LEVEL DENSITY OF ELLIPTIC CURVE $L$-FUNCTIONS }}
\vspace {2 in}
\author{D.\ K.\ Huynh, J.\ P.\ Keating and N.\ C.\ Snaith\\
        School of Mathematics,\\ University of Bristol,\\
            Bristol BS8 1TW, UK}
\date{\today}
\maketitle \thispagestyle{empty} \vspace{.5cm}
\begin{abstract}
It is believed that, in the limit as the conductor tends to infinity, correlations between the zeros of elliptic curve $L$-functions averaged within families follow the distribution laws of the eigenvalues of random matrices drawn from
the orthogonal group.  For test functions with restricted support, this is known to be the true for
the one- and two-level densities of
zeros within the families studied to
date. However, for finite conductor Miller's
experimental data reveal an interesting discrepancy from
these limiting results. Here we use the $L$-functions ratios
conjectures  to calculate the 1-level density for the family of
even quadratic twists of an elliptic curve $L$-function for large but finite conductor.  This
gives a formula for the leading and lower order terms up to an error term that is conjectured to be significantly smaller. The lower order terms explain many of
the features of the zero statistics for relatively small conductor and
model the very slow convergence to the infinite conductor limit.
 However, our main observation is that they do not
 capture the behaviour of zeros in the important region
 very close to the critical point and so do not explain Miller's discrepancy.   This therefore implies
 that a more accurate model for statistics near to this point
 needs to be developed.
\end{abstract}


\section{ Introduction}

The conjecture that the limiting statistical properties of the zeros of $L$-functions
 may be modeled by those of the eigenvalues of random matrices goes back to
 Montgomery \cite{kn:mont73}, who introduced it in the context of the Riemann
 zeta-function.  For the Riemann zeros this conjecture is supported by extensive
 numerical  \cite{kn:odlyzko97} and theoretical
 \cite{kn:mont73,kn:hejhal94,kn:bogkea95,kn:bogkea96,kn:rudsar} calculations.
 The generalization to zero statistics within families of $L$-functions
 was developed by Katz and Sarnak \cite{kn:katzsarnak99a,kn:katzsarnak99b},
 and again there is much evidence supporting it \cite{kn:rub01}.
 Random matrix models for the moments of the Riemann zeta-function on its
 critical line and for central values of $L$-functions within families were
 introduced by Keating and Snaith \cite{kn:keasna00a,kn:keasna00b}, and have
 since been developed extensively
 \cite{kn:confar00,kn:cfkrs,kn:ghk,kn:buikea07,kn:buikea08,kn:cfkrs2}.
 For more background, see \cite{kn:mezzsna}.

The random-matrix moment conjectures extend naturally to ratios of $L$-functions.
The $L$-functions ratios conjectures were stimulated by the work of
Farmer, who, in 1995, made a conjecture for shifted moments of the
Riemann zeta-function \cite{kn:farmer95}. Nonnenmacher and
Zirnbauer \cite{kn:nonzir02} found formulas for the ratios of
characteristic polynomials of random matrices coming from one of
the classical compact groups. This was formalised and written up
by Conrey, Farmer and Zirnbauer \cite{kn:cfz1} and lead to the
development of corresponding ratios conjectures for $L$-functions
in number theory \cite{kn:cfz2}.

The Birch/Swinnerton-Dyer conjecture asserts that the rank of an
elliptic curve is equal to the order of vanishing at the central
point of the associated $L$-function.  The idea of using random
matrix theory to predict the frequency of non-zero rank in
families of elliptic curves was introduced by Conrey, Keating,
Rubinstein and Snaith \cite{kn:ckrs00,kn:ckrs05}. An interesting
extension of this is to find a random matrix model for elliptic
curve $L$-functions of
 a given order of vanishing at the critical point.  The first steps in
 this direction have been taken by Snaith \cite{kn:snaith05a} and
 Miller/Due\~{n}ez \cite{kn:mil05}, but it is clear from Miller's numerical computations
 that there is a still simpler problem concerning the zero statistics of families of
 rank zero curves that is far from being
 understood.  This problem is the main motivation for the work we shall report on here.

According to the Katz/Sarnak philosophy \cite{kn:katzsarnak99a,kn:katzsarnak99b},
zeros of families of $L$-functions show the same statistical
behaviour as eigenvalues of random matrices drawn from one of the
classical compact groups. The zeros of a family
of elliptic curve $L$-functions with even (odd) functional
equation should follow the distribution laws of eigenvalues of the even
(odd) orthogonal group.  Rigorous calculations \cite{kn:miller02,
kn:mil04,kn:young}  show that as the conductor (the parameter that
orders $L$-functions within a family) tends to infinity, the one-
and two-level densities do indeed tend to the expected orthogonal
forms for several different families of elliptic curves.
That is, as the conductor tends to infinity, the zero statistics
approach the scaling limit for large matrix size of the
corresponding statistic for the eigenvalues of matrices from $SO(2N)$
or $SO(2N+1)$. (Similar agreement with random matrix theory is shown for many
other families of $L$-functions, see for example
\cite{kn:duemil06,
kn:fouiwa03,kn:guloglu05,kn:hugrud03a,kn:hugmil07,kn:ILS99,kn:ozlsny99,
kn:ricroy08a,kn:royer01,kn:rub01}.)  The test functions involved in these
calculations have a limited range of support, but nonetheless the
evidence is compelling. Thus it was surprising to see in
Miller's numerical results \cite{kn:mil05} a distinct repulsion of the zeros from the
central point for a family of $L$-functions of rank 0
elliptic curves, because no repulsion is seen in the statistics of $SO(2N)$
eigenvalues.   Of course, in numerical computations the
conductor is finite, and so it is clear that an explanation is needed for finite conductor statistics and how they approach
 the limiting
$SO(2N)$ statistic.

We do have a relatively complete understanding of the way in which the random matrix limit is approached for the zero statistics of the Riemann zeta function at a height $T$ up the critical line as $T\rightarrow\infty$.  Berry first wrote down an approximate formula describing the finite-$T$ corrections to the random matrix limiting form for a statistic related to the 2-point correlation function in \cite{kn:berry88} and showed that this described Odlyzko's data remarkably accurately.  Later, a formula that is believed to capture all of the essential features was derived by Bogomolny and Keating \cite{kn:bk96}.  The terms in the Bogomolny-Keating formula that describe the corrections to the random matrix limit are often referred to as {\it lower order terms}.  See \cite{kn:berrykeating99} for an overview and numerical illustrations.  More recently, Conrey and Snaith \cite{kn:consna06, kn:consna08} have shown how the Bogomolny-Keating formula and its extension to all $n$-point correlation functions can be recovered from the $L$-functions ratios conjectures \cite{kn:cfz2}.   There have also been
investigations of lower order terms in the zero statistics of
various families of $L$-functions
\cite{kn:fouiwa03,kn:mil07,kn:mil08,kn:ricroy08b,kn:young05}.  In particular, Conrey and Snaith have shown how such terms can also be recovered from the ratios conjectures \cite{kn:consna06}.  It is thus natural in this context to seek the explanation for the surprising discrepancy observed by Miller in these lower order terms.

\begin{figure}[h]
\begin{center}
\hspace*{65mm}\includegraphics[scale=0.8]{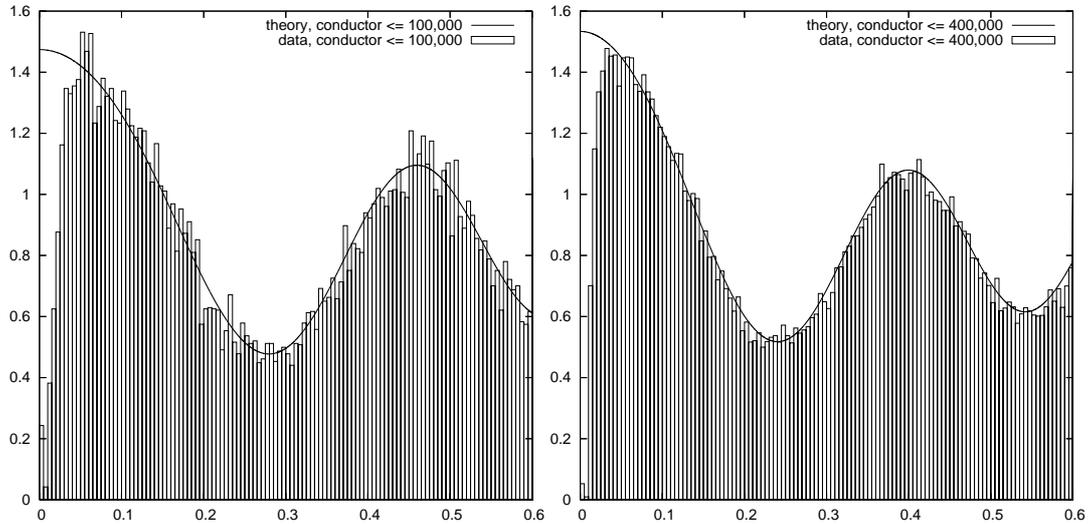}
\caption{1-level density of unscaled zeros from 0 up to height 0.6
of even quadratic twists of $L_{E_{11}}$ with $0 < d < 100,000$
for {\it left} and $0 < d < 400,000$ for {\it right} hand side,
prediction (solid), from (\ref{oneleveldensity}), versus numerical
data (bar chart).
} \label{fig:nrzero}
\end{center}
\end{figure}

In this paper we examine lower order terms in the 1-level
density of the zeros of a family of elliptic curve $L$-functions.
Specifically, we investigate even quadratic twists of an
elliptic curve $L$-function, for which we calculate the zeros
numerically with Rubinstein's {\sffamily lcalc}
\cite{kn:rubinstein}. Using the ratios
conjectures we derive a formula for the 1-level density that describes convincingly the
 intricate structure of the numerical data away from the central point and so explains the rate of approach to the
 random matrix limit in this region.  However, most interestingly, our formula fails to describe the
 region very close to the central point.  To illustrate our main results, we plot in
 figure \ref{fig:nrzero} a numerical evaluation of the 1-level density together with our formula.  Miller's discrepancy corresponds to the region near to the origin.  Our main conclusion here is then that the
 explanation for the zero distribution in this region lies beyond the models combining random matrix theory and arithmetical lower order terms considered so far; that is, these formulae are not sufficient to explain the discrepancy.  We plan to explore augmented models that build on the present calculation to explain the phenomenon in a future paper with E. Due{\~n}ez
and S. J. Miller


\section{ The 1-level density formula}

Let the $L$-function $L_E(s)$ associated with an elliptic curve
$E$ be given by the Dirichlet series
\begin{equation}
L_E(s)=\sum_{n=1}^{\infty} \frac{\lambda(n)}{n^s},
\end{equation}
where the coefficients ($\lambda(n)=a(n)/\sqrt{n}$, with $a_p=p+1-\#E({\mathbb F}_p)$,
$\#E({\mathbb F}_p)$ being the number of points on $E$ counted over
${\mathbb F}_p$) have been normalised
 so that the functional equation relates $s$ to
$1-s$:
\begin{equation}
 L_E(s)=\omega(E) \left(\frac{2\pi}{\sqrt{M}}\right)^{2s-1}
\frac{\Gamma(3/2-s)}{ \Gamma(s+1/2)}L_E(1-s).
\end{equation}
Here $M$ is the conductor of the elliptic curve $E$; we will
consider only prime $M$. Also, $\omega(E)$ is $+1$ or $-1$
resulting, respectively, in an even or odd functional equation for
$L_E$.

 Let $L_E(s, \chi_d)$ denote the $L$-function obtained by twisting  $L_E(s)$ quadratically.
Here $d$ is a fundamental discriminant, i.e., $d \in
\mathbb{Z}-\{1\}$, s.t. $p^2 \nmid d$ for all odd primes $p$ and
$d \equiv 1 \modd 4$ or $d \equiv 8, 12 \modd 16$, and $\chi_d$ is
the Kronecker symbol. Then the twisted $L$-function (which is
itself the $L$-function associated with another elliptic curve
$E_d$) is given by
\begin{equation} \label{twistedLfunction}
L_E(s, \chi_d) = \sum_{n = 1}^\infty \frac{\lambda(n)
\chi_d(n)}{n^s} = \prod_p \left(1 -
\frac{\lambda(p)\chi_d(p)}{p^s} +
\frac{\psi_{M}(p)\chi_d(p)^2}{p^{2s}} \right)^{-1}
\end{equation}
where $\psi_{M}$ is the principal Dirichlet character of modulus
$M$:
\begin{equation}
\psi_{M}(p) =
\begin{cases}
1 \text{~if~} p \nmid M\\
0 \text{~otherwise}.
\end{cases}
\end{equation}
The functional equation of this $L$-function is
\begin{equation}
L_E(s,\chi_d)=\chi_d(-M)\omega(E)
\left(\frac{2\pi}{\sqrt{M}|d|}\right)^{2s-1} \frac{\Gamma(3/2-s)}{
\Gamma(s+1/2)}L_E(1-s,\chi_d).
\end{equation}

In order to derive the 1-level density of the zeros near the
critical point $s=1/2$ of $L$-functions in this family of
quadratic twists, we consider the average over the family of a
ratio of $L$-functions evaluated at different points:
\begin{equation}
R_E(\alpha, \gamma) := \label{RE} \sum_{\substack{0<d \leq X\\
\chi_d(-M)\omega_E=+1}} \frac{L_E(1/2 + \alpha, \chi_d)}{L_E(1/2 +
\gamma, \chi_d)}.
\end{equation}
This is an average over those twisted $L$-functions that have even
functional equations and $0< d \leq X$.  Requiring an even
functional equation imposes a restriction on $d \mod M$.  We
follow the recipe of \cite{kn:cfkrs}, \cite{kn:cfz2} and the
calculations in \cite{kn:consna06} to derive a formula for
$R_E(\alpha, \gamma)$ via the ratios conjecture.  Note that
arriving at a ratios conjecture  entails applying a list of
manipulations, several of which introduce errors large enough to
be significant.  The miracle is that these errors appear to cancel
out and the recipe yields formulae that have been checked
numerically and against specific known cases in many different
situations (see \cite{kn:cfz2,kn:consna06}). Recent work of Steven
J. Miller \cite{kn:mil07} has shown that a rigorous calculation of
the 1-level density for the family of real quadratic Dirichlet
$L$-functions matches exactly, for a suitably chosen test
function, the prediction obtained by applying the ratios recipe.
See also \cite{kn:sto08} for further investigations of the ratios
conjecture and the 1-level density of the same family of Dirichlet
$L$-functions and \cite{kn:mil08} for Miller's extension of
\cite{kn:mil07} to families of cuspidal newforms.

We use (\ref{twistedLfunction}) to replace $L_E(s, \chi_d)$ in the
denominator of (\ref{RE}) by
\begin{equation}
\frac{1}{L_E(s, \chi_d)} = \sum_{n = 1}^\infty \frac{\mu_E(n)\chi_d(n)}{n^s}
\end{equation}
where $\mu_E(n)$ is a multiplicative function defined as
\begin{equation}
\mu_E(n) =
\begin{cases}
-\lambda(p), \mbox{if~} n = p\\
\psi_{M}(p), \mbox{if~} n = p^2\\
0, \mbox{if~} n = p^k, k > 2.\\
\end{cases}
\end{equation}
We use the approximate functional equation for the $L$-function in
the numerator of (\ref{RE}):
\begin{align} \nonumber
L_E(1/2 + \alpha, \chi_d) = & \sum_{m < x} \frac{\chi_d(m)
\lambda(m)}{m^{1/2 + \alpha}} +
\left(\frac{\sqrt{M}|d|}{2\pi}\right)^{-2\alpha} \frac{\Gamma(1-
\alpha)}{\Gamma(1 + \alpha)} \sum_{n < y} \frac{\chi_d(n)
\lambda(n)}{n^{1/2-\alpha}} \\& + \mbox{remainder},
\label{approximatefunctionalequation}
\end{align}
where $M$ is the conductor of the elliptic curve $E$ and $xy = d^2 / (2\pi)$. Therefore using the first sum
of the approximate functional equation (\ref{approximatefunctionalequation}) we get
\begin{equation} \label{RE1}
R_E^1(\alpha, \gamma) := \sum_{\substack{0<d \leq X\\
\chi_d(-M)\omega_E=+1}} \sum_{h, m} \frac{\lambda(m) \mu_E(h)
\chi_d(mh)}{m^{1/2 + \alpha} h^{1/2 + \gamma}}.
\end{equation}
We denote by $R_E^2(\alpha, \gamma)$ the expression that results
from using the second sum in the approximate functional equation
(\ref{approximatefunctionalequation}). Thus
\begin{equation} \label{RE_sum}
R_E(\alpha, \gamma) \approx R_E^1(\alpha,\gamma) + R_E^2(\alpha, \gamma).
\end{equation}
The ratios recipe now calls for a replacement of $\chi_d(mh)$ with
its average over the family (the set of $d$'s being summed over).
We set
\begin{equation}
X^*= \sum_{\substack{0<d \leq X\\
\chi_d(-M)\omega_E=+1}}1 {\rm\;\;\;and\;\;\;} X_b^*=
\sum_{\substack{0<d \leq X\\d=b\mod M}}1 \end{equation}
 as the number of
fundamental discriminants below $X$ that we are summing over and
note (see \cite{kn:cfkrs}, Theorem 3.1.1)
\begin{equation}\label{eq:harmonic}
\frac{1}{X^*_b}\sum_{\substack{0<d \leq X \\d=b\mod M}} \chi_d(n)
\approx\left\{
\begin{array}{cl}
\chi_b(g)a(n) & \mbox{if }n=g\square, \mbox{ with }
(\square,M)=1\mbox{ and if all prime}\\&  \mbox{factors of $g$ are
prime factors of $M$} \\0& \mbox{otherwise,}
\end{array}\right.
\end{equation}
where
\begin{equation}
a(n) = \prod_{p|\square} \frac{p}{p+1}.
\end{equation}
 This is to say that terms not of the form $n=g\square$ can be disregarded (this is the so-called `harmonic
detector' which is mentioned in \cite{kn:cfkrs}).  Since we are
considering only curves with prime conductor $M$, $g$ is simply a
power of $M$.  Note that in the cases we are interested in
$\chi_b(g)=\omega_E^{\ell}$ for $g=M^{\ell}$ because $d$ has been
chosen such that $\chi_b(M)=\chi_d(M)=\omega_E$ (we have
$\chi_d(M)=\chi_d(-M)$ since we are considering only positive
$d$).

Concentrating on $R_E^1$, we replace $\chi_{d}(mh)$ with the
average given by (\ref{eq:harmonic}) and so restrict the sum as
follows:
\begin{equation}
R_E^1(\alpha, \gamma) \approx X^* \sum_{hm = \square M^{\ell}}
\frac{\lambda(m) \mu_E(h) a(mh)\omega_E^{\ell}}{m^{1/2 +
\alpha}h^{1/2 + \gamma}},
\end{equation}
with $(\square,M)=1$ and $g$ divisible only by primes dividing
$M$.  We write this sum as an Euler product (for convenience
denoting by $h$ the exponent on primes dividing $h$ in the sum
above and similarly for $m$) and note that if $m+h\geq 1$ then
$a(p^{m+h}) = p/p+1$ for primes not dividing the conductor,
whereas $a(p^{m+h})=1$ if the prime does divide the conductor. So
we obtain
\begin{equation}
R_E^1(\alpha, \gamma) \approx X^* V_|(\alpha,
\gamma)V_{\nmid}(\alpha,\gamma)
\end{equation}
where
\begin{eqnarray}\label{summe}
&&V_\nmid(\alpha, \gamma) := \prod_{p\nmid M} \Bigg(1 +
\frac{p}{p+1} \sum_{\substack{m,h\geq 0\\m + h > 0 \\ m + h ~{\rm
even}}} \frac{\lambda(p^m) \mu_E(p^h)}{p^{m(1/2 + \alpha) + h(1/2
+
\gamma)}} \Bigg)\\
&&V_{|}(\alpha,\gamma):= \prod_{p| M}\Bigg(\sum_{h,m\geq 0}
\frac{\lambda(p^m)\mu_E(p^h)\omega_E^{m+h}} {p^{m(1/2 + \alpha) +
h(1/2 + \gamma)}}\Bigg).\label{eq:vdivide}
\end{eqnarray}
Since $\mu_E(p^h) = 0$ for most powers of $p$, we only need to
consider $h = 0,1,2$ in the sum in (\ref{summe}) and $h=0,1$ in
(\ref{eq:vdivide}). Then the Euler products become
\begin{eqnarray} \label{Euler}
V_\nmid(\alpha, \gamma) &= &\prod_{p\nmid M} \left(1 + \frac{p}{p
+ 1}\left(\sum_{m = 1}^\infty \frac{\lambda(p^{2m})}{p^{m(1 +
2\alpha)}} - \frac{\lambda(p)}{p^{1 + \alpha + \gamma}}\sum_{m =
0}^\infty \frac{\lambda(p^{2m+1})}{p^{m(1 +
2\alpha)}}\right.\right.
\nonumber\\
&&\qquad\qquad\qquad\left.\left.+ \frac{1}{p^{1 + 2\gamma}}\sum_{m
= 0}^\infty \frac{\lambda(p^{2m})}{p^{m(1 + 2\alpha)}} \right)
\right)
\end{eqnarray}
and
\begin{eqnarray}\label{eq:Eulerdivide}
V_{|}(\alpha,\gamma)=\prod_{p| M}\Bigg(\sum_{m= 0}^{\infty}\bigg(
\frac{\lambda(p^m)\omega_E^{m}} {p^{m(1/2 + \alpha)
}}-\frac{\lambda(p)\lambda(p^m)\omega_E^{m+1}}{p^{m(1/2+\alpha)+1/2+\gamma}}\bigg)\Bigg).
\end{eqnarray}
We now factor out the divergent part of $R_E^1$ using the Riemann
zeta function and also, for convenience, we will factor out the
symmetric square $L$-function associated with $L_E$.  This leaves
us with a convergent Euler product. In the following, for
simplicity, we shall only deal with elliptic curves with prime
conductor, $M$. Recall that the Euler product of a Hasse-Weil
$L$-function $L_E(s)$ coming from the elliptic curve $E$, with
Dirichlet coefficients $\lambda(n)$ normalised so that the
functional equation relates $s$  to $1 - s$, has the form
\begin{equation}
L_E(s) = \prod_{p|M} (1 - \lambda(p) p^{-s})^{-1} \prod_{p\nmid M}
(1 - \lambda(p) p^{-s} + p^{-2s})^{-1}.
\end{equation}
 Now we can write this product as
\begin{equation}
L_E(s) = \prod_p (1 - \alpha(p)p^{-s})^{-1} (1 - \beta(p)p^{-s})^{-1}
\end{equation}
where
\begin{equation} \label{alphaplusbeta}
\alpha(p) + \beta(p) = \lambda(p)
\end{equation}
and
\begin{equation}
\alpha(p) \beta(p) =
\begin{cases}
0 \mbox{~for~} p | M\\
1 \mbox{~for~} p \nmid M. \label{Steve}
\end{cases}
\end{equation}
Let $L_E(\sym^2, s)$ denote the symmetric square $L$-function.
Then by definition (see \cite{kn:iwaniec}, page 251)
\begin{equation} \label{Iwaniec}
L_E(\sym^2, s) = \prod_p (1 - \alpha^2(p)p^{-s})^{-1} (1 - \alpha(p)\beta(p)p^{-s})^{-1}
(1 - \beta^2(p)p^{-s})^{-1}.
\end{equation}
We have (see \cite{kn:conrey04b}, page 236)
\begin{equation}
\lambda(m) \lambda(n) = \sum_{\substack{d|(m, n)\\(d, M) = 1}} \lambda(mn/d^2),
\end{equation}
(where $M$ is the conductor of $E$)  and in particular we have for
$p \nmid M$
\begin{eqnarray} \label{quadrat}
\lambda(p)^2 & = & \lambda(p^2) + 1\\ \label{lang}
\lambda(p^{2m +1})\lambda(p) & = & \lambda(p^{2m +2}) + \lambda(p^{2m}).
\end{eqnarray}
We wish to write the Euler product in (\ref{Iwaniec}) in terms of $\lambda(p)$, so we start
by using (\ref{alphaplusbeta}) to obtain
\begin{eqnarray}
L_E(\sym^2, s)& =& \prod_p \left(1 - \frac{\lambda(p)^2 -
\alpha(p)\beta(p)}{p^s} + \frac{\alpha(p)\beta(p)(\lambda(p)^2 -
\alpha(p)\beta(p))}{p^{2s}} \right.\nonumber \\
&&\qquad\qquad\qquad\left.- \frac{(\alpha(p)\beta(p))^3}{p^{3s}}
\right)^{-1}.
\end{eqnarray}
We now distinguish between $p|M$ and $p \nmid M$, and so, using
(\ref{quadrat}) and (\ref{Steve}), we have
\begin{equation} \label{symmetric_old}
L_E(\sym^2, s) = \prod_{p | M} \left(1 -
\frac{\lambda(p)^2}{p^s}\right)^{-1} \prod_{p \nmid M} \left(1 -
\frac{\lambda(p^2)}{p^s} + \frac{\lambda(p^2)}{p^{2s}} -
\frac{1}{p^{3s}} \right)^{-1}.
\end{equation}

\kommentar{
Remark: Note that for $p | \Delta$ we have $-p^{1/2}\lambda(p) = \pm 1 \imply p\lambda(p)^2 = 1 \imply \lambda(p)^2 = 1/p$.
So we can write (\ref{symmetric}) also as
$$
L_E(\sym^2, s) =
\prod_{p | \Delta} \left(1 - \frac{1}{p^{s+1}}\right)^{-1}
\prod_{p \nmid \Delta} \left(1 - \frac{\lambda(p^2)}{p^s} + \frac{\lambda(p^2)}{p^{2s}} - \frac{1}{p^{3s}} \right)^{-1}.
$$
}

Now we reconsider the Euler products in (\ref{Euler}) and
(\ref{eq:Eulerdivide}). In constructing ratios conjectures we
usually allow $-\tfrac{1}{4}<{\rm Re} \alpha<\tfrac{1}{4}$ and
$\log X\ll {\rm Re}\gamma<\tfrac{1}{4}$, where the bounds at
$\tfrac{\pm 1}{4}$ allow us to control the convergence of Euler
products of the type (\ref{Euler}). In fact, in this application
the real parts of $\alpha$ and $\gamma$ can be considered as very
small. Thus we can write
\begin{eqnarray}
V_{\nmid}(\alpha, \gamma) &=& \prod_{p\nmid M} \left(1 +
\frac{p}{p + 1}\left(\sum_{m = 1}^\infty
\frac{\lambda(p^{2m})}{p^{m(1 + 2\alpha)}} -
\frac{\lambda(p)}{p^{1 + \alpha + \gamma}}\sum_{m = 0}^\infty
\frac{\lambda(p^{2m+1})}{p^{m(1 + 2\alpha)}}
\right.\right.\nonumber\\
&&\qquad\qquad\qquad\qquad \left.\left.+ \frac{1}{p^{1 +
2\gamma}}\sum_{m = 0}^\infty \frac{\lambda(p^{2m})}{p^{m(1 +
2\alpha)}} \right)
\right)\nonumber \\
&=& \label{orderforeulerproduct} \prod_{p\nmid M} \left(1 +
\frac{\lambda(p^2)}{p^{1 + 2\alpha}} - \frac{\lambda(p^2) + 1
}{p^{1 + \alpha + \gamma}} + \frac{1}{p^{1 + 2 \gamma}} +
\cdots\right),
\end{eqnarray}
where the $\cdots$ indicate terms that converge like $1/p^2$ when
$\alpha$ and $\gamma$ are small.  We now use the following
approximations to factor out the divergent or slowly converging
terms. By (\ref{symmetric_old}) we have
\begin{equation}
L_E(\sym^2, 1 + 2\alpha)=\prod_p\left( 1+\frac{\lambda(p^2)}{p^{1
+ 2\alpha}} +\cdots \right)
\end{equation}
and
\begin{equation}
\frac{1}{L_E(\sym^2, 1+\alpha + \gamma)} \frac{1}{\zeta(1 + \alpha
+ \gamma)}=\prod_p\left(1- \frac{\lambda(p^2) + 1 }{p^{1 + \alpha
+ \gamma}}+\cdots\right).
\end{equation}
Also, since there is only one prime that divides the conductor
$M$, a factor of $\zeta(1+2\gamma)$ will account for the
divergence of the term $\frac{1}{p^{1 + 2 \gamma}}$ in
(\ref{orderforeulerproduct}).

Hence we can write
\begin{equation}
V_\nmid(\alpha, \gamma)V_|(\alpha, \gamma)=Y_E(\alpha,
\gamma)A_E(\alpha, \gamma),
\end{equation}
where \begin{equation}\label{eq:YE}
 Y_E(\alpha, \gamma) = \frac{\zeta(1 +
2\gamma) L_E(\sym^2, 1 + 2\alpha)}{\zeta(1 + \alpha + \gamma)
L_E(\sym^2, 1+ \alpha+\gamma)}.
\end{equation}
$A_E(\alpha, \gamma)$ is given by
\begin{align} \nonumber
&A_E(\alpha, \gamma) = ~ Y_E^{-1}(\alpha, \gamma)\times
\prod_{p\nmid M} \left(1 + \frac{p}{p + 1}\left(\sum_{m =
1}^\infty \frac{\lambda(p^{2m})}{p^{m(1 + 2\alpha)}}\right.\right.
\\ \label{AE} &\left.\left. \qquad\qquad- \frac{\lambda(p)}{p^{1 + \alpha + \gamma}}\sum_{m =
0}^\infty \frac{\lambda(p^{2m+1})}{p^{m(1 + 2\alpha)}} +
\frac{1}{p^{1 + 2\gamma}}\sum_{m = 0}^\infty
\frac{\lambda(p^{2m})}{p^{m(1 + 2\alpha)}} \right) \right)\\
\nonumber &\qquad\qquad\qquad\qquad\times\prod_{p|
M}\Bigg(\sum_{m= 0}^{\infty}\bigg( \frac{\lambda(p^m)\omega_E^{m}}
{p^{m(1/2 + \alpha)
}}-\frac{\lambda(p)}{p^{1/2+\gamma}}\frac{\lambda(p^m)\omega_E^{m+1}}{p^{m(1/2+\alpha)}}\bigg)\Bigg)
\end{align}
and is analytic as $\alpha, \gamma \rightarrow 0$. Hence, by
recalling (\ref{RE1}), we find
\begin{equation}
R_E^1(\alpha,
\gamma) \approx \sum_{\substack{0<d \leq X\\
\chi_d(-M)\omega_E=+1}} Y_E(\alpha, \gamma) A_E(\alpha,
\gamma).
\end{equation}
We obtain the other sum $R_E^2(\alpha, \gamma)$ in (\ref{RE_sum})
by using the second term in the approximate functional equation
(\ref{approximatefunctionalequation}) and carrying out exactly the
same steps as above:
\begin{equation}
R_E^2(\alpha, \gamma)\approx \sum_{\substack{0<d \leq X\\
\chi_d(-M)\omega_E=+1}} \left(\frac{\sqrt{M}|d|}{2\pi}
\right)^{-2\alpha} \frac{\Gamma(1 - \alpha)}{\Gamma(1 + \alpha)}
Y_E(-\alpha, \gamma) A_E(-\alpha, \gamma).
\end{equation}

By applying the ratios conjecture recipe, we therefore have the
result:
\begin{conjecture}[Ratios Conjecture]\label{conj:ratios}  For some reasonable conditions such as $-\frac{1}{4}<\RRe
\alpha<\frac{1}{4}$, $\frac{1}{\log X} \ll \RRe
\gamma<\frac{1}{4}$ and $\IIm \alpha,\IIm\gamma\ll
X^{1-\varepsilon}$, we have
\begin{align}
R_E(\alpha, \gamma) = &\sum_{\substack{0<d \leq X\nonumber\\
\chi_d(-M)\omega_E=+1}} \frac{L_E(1/2 + \alpha, \chi_d)}{L_E(1/2 +
\gamma, \chi_d)}
 \\
= & \sum_{\substack{0<d \leq X\nonumber\\
\chi_d(-M)\omega_E=+1}} \left( Y_E A_E(\alpha, \gamma) + \left(\frac{\sqrt{M}|d|}{2\pi}\right)^{-2\alpha} \frac{\Gamma(1 - \alpha)}{\Gamma(1 + \alpha)}Y_E A_E(-\alpha, \gamma)\right) \\
& \qquad\qquad\qquad+ O(X^{1/2 + \varepsilon}),\nonumber
\end{align}
where $Y_E$ and $A_E$ are defined at (\ref{eq:YE}) and (\ref{AE}),
respectively, $M$ is the (prime) conductor of the $L$-function
$L_E(s)$ and $\omega_E$ is the sign from its functional equation.
\end{conjecture}

We note that the error term $O(X^{1/2+\varepsilon})$ is part of
the statement of the ratios conjecture; the power on $X$ is not
suggested by any of the steps used in arriving at the main
expression in Conjecture \ref{conj:ratios}.  At the end of Section
\ref{sect:data} we propose that the limited data we have available
supports a power saving on the error term, but not necessarily a
power of 1/2.

To calculate the 1-level density we actually need the average of
the logarithmic derivative of $L$-functions in this family, so we
note that
\begin{eqnarray}
\sum_{\substack{0<d \leq X\\
\chi_d(-M)\omega_E=+1}} \frac{L_E'(1/2 + r, \chi_d)}{L_E(1/2 + r,
\chi_d)} = \frac{d}{d \alpha}R_E(\alpha, \gamma)\Big|_{\alpha =
\gamma = r}.
\end{eqnarray}
Using (\ref{lang}) for primes not dividing $M$ and the
multiplicativity of $\lambda(p)$ for $p|M$, we get $A_E(r,r) = 1$
and we have, with
\begin{equation}
\label{eq:AE1}
A_E^1(r,r)=\frac{d}{d\alpha}A_E(\alpha,\gamma)\big|_{\alpha=\gamma=r},
\end{equation}
\begin{align}\nonumber
\frac{d}{d \alpha} Y_E A_E(\alpha, \gamma)\Big|_{\alpha = \gamma = r}
= & -\frac{\zeta'(1 + 2r)}{\zeta(1 + 2r)} A_E(r,r) + \frac{L_E'(\sym^2, 1 + 2r)}{L_E(\sym^2, 1+2r)}A_E(r,r) + A_E^1(r,r) \\
= & -\frac{\zeta'(1 + 2r)}{\zeta(1 + 2r)} + \frac{L_E'(\sym^2, 1 +
2r)}{L_E(\sym^2, 1+2r)} + A_E^1(r,r)
\end{align}
and
\begin{align}\nonumber
\frac{d}{d \alpha}
\left(\frac{\sqrt{M}|d|}{2\pi}\right)^{-2\alpha} \frac{\Gamma(1 -
\alpha)}{\Gamma(1 + \alpha)}\{Y_E(-\alpha, \gamma) A_E(-\alpha,
\gamma)\}\Big|_{\alpha = \gamma = r}
\\= - \left(\frac{\sqrt{M}|d|}{2\pi}\right)^{-2r} \frac{\Gamma(1 - \alpha)}{\Gamma(1 + \alpha)}
\frac{\zeta(1 + 2r)L_E(\sym^2, 1 - 2r)}{L_E(\sym^2,1)}A_E(-r, r).
\end{align}
Therefore we have for the logarithmic derivative the following:
\begin{theorem}\label{theo:logderiv} Assuming the Ratios
Conjecture \ref{conj:ratios} and $\frac{1}{\log X} \ll \RRe
($r$)<\frac{1}{4}$ and $\IIm ($r$)\ll X^{1-\varepsilon}$, the
average of the logarithmic derivative over a family of quadratic
twists (with even functional equation) of the $L$-function of an
elliptic curve with prime conductor $M$ is
\begin{align} \nonumber
& \sum_{\substack{0<d \leq X\\
\chi_d(-M)\omega_E=+1}} \frac{L_E'(1/2 + r, \chi_d)}{L_E(1/2 + r,
\chi_d)}
\\  & = \sum_{\substack{0<d \leq X\\
\chi_d(-M)\omega_E=+1}} \Bigg( -\frac{\zeta'(1 + 2r)}{\zeta(1 +
2r)} + \frac{L_E'({\rm sym}^2, 1 + 2r)}{L_E({\rm sym}^2, 1+2r)} +
A_E^1(r,r)
\\ \label{eins} & ~~~~-
\left(\frac{\sqrt{M}|d|}{2\pi}\right)^{-2r} \frac{\Gamma(1 -
r)}{\Gamma(1 + r)} \frac{\zeta(1 + 2r)L_E({\rm sym}^2, 1 -
2r)}{L_E({\rm sym}^2,1)}A_E(-r, r) \Bigg) +
O(X^{1/2+\varepsilon}).\nonumber
\end{align}
Here $\omega_E$ is the sign from the functional equation of $L_E$,
$L_E({\rm sym}^2,s)$ is the associated symmetric square
$L$-function (defined at (\ref{Iwaniec})), and $A_E$ and $A_E^1$
are arithmetic factors defined at (\ref{AE}) and (\ref{eq:AE1}),
respectively.
\end{theorem}

Let $\gamma_d$ denote the ordinate of a generic zero of $L_E(s, \chi_d)$ on the half line.
We consider the 1-level density
\begin{equation}
S_1(f) := \sum_{\substack{0<d \leq X\\
\chi_d(-M)\omega_E=+1}} \sum_{\gamma_d} f(\gamma_d)
\end{equation}
where $f$ is some nice test function, say an even Schwartz function.
By the argument principle we have
\begin{equation}
S_1(f) = \sum_{\substack{0<d \leq X\\
\chi_d(-M)\omega_E=+1}} \frac{1}{2\pi i} \left( \int_{(c)} -
\int_{(1-c)}\right) \frac{L'(s, \chi_d)}{L(s, \chi_d)} f (-i(s -
1/2))ds
\end{equation}
where $(c)$ denotes a vertical line from $c -i\infty$ to $c + i\infty$ and $3/4 > c > 1/2 + 1/\log X$.
The integral on the $c$-line is
\begin{equation}
\frac{1}{2\pi} \int_{-\infty}^\infty f(t - i(c-1/2)) \sum_{\substack{0<d \leq X\\
\chi_d(-M)\omega_E=+1}} \frac{L_E'(1/2 + (c -1/2
+it),\chi_d)}{L_E(1/2 + (c-1/2 +it), \chi_d)}dt.
\end{equation}
The sum over $d$ can be replaced by Theorem \ref{theo:logderiv}.
The bounds on the size $t$ coming from the ratios conjecture
should not limit us here.  It is not entirely known in what range
of the parameters the ratios conjecture holds, but the test
function $f$ can be chosen to decay sufficiently fast that the
tails of the integrand, where the ratios conjecture might fail,
will not contribute significantly. (See the 1-level density
section of \cite{kn:consna06} for more detailed analysis.) Next we
move the path of integration to $c = 1/2$ as the integrand is
regular at $t = 0$  and get
\begin{align}
& \frac{1}{2\pi} \int_{-\infty}^\infty f(t) \sum_{\substack{0<d \leq X\\
\chi_d(-M)\omega_E=+1}} \Bigg( -\frac{\zeta'(1 + 2it)}{\zeta(1 +
2it)} + \frac{L_E'(\sym^2, 1 + 2it)}{L_E(\sym^2, 1+2it)} +
A_E^1(it,it)
\\ & \qquad - \left(\frac{\sqrt{M}|d|}{2\pi}\right)^{-2it} \frac{\Gamma(1 - it)}{\Gamma(1 + it)}
\frac{\zeta(1 + 2it)L_E(\sym^2, 1 - 2it)}{L_E(\sym^2,1)}A_E(-it,
it) \Bigg)dt \nonumber \\
&\qquad\qquad\qquad\qquad\qquad
\qquad\qquad\qquad\qquad\qquad\nonumber+ O(X^{1/2+\varepsilon}).
\end{align}

For the integral on the line with real part $1-c$, we use the
functional equation
\begin{equation}
L_E(s, \chi_d) = \chi_d(-M)\omega_E X(s, \chi_d) L_E(1-s, \chi_d)
\end{equation}
with
\begin{equation} \label{xi_equation}
X(s, \chi_d) =
\left(\frac{\sqrt{M}|d|}{2\pi}\right)^{1-2s} \frac{\Gamma(3/2-s)}{
\Gamma(s+1/2)}
\end{equation}
to obtain
\begin{equation}
\frac{L_E'(1-s, \chi_d)}{L_E(1-s, \chi_d)} = \frac{X'(s, \chi_d)}{X(s, \chi_d)} - \frac{L_E'(s, \chi_d)}{L_E(s, \chi_d)}. \label{zwei}
\end{equation}
The logarithmic derivative of (\ref{xi_equation}) evaluated
at $s = 1/2 + \alpha$ is
\begin{equation}
\frac{X'(1/2 + \alpha, \chi_d)}{X(1/2 + \alpha, \chi_d)} = -2\log
\left(\frac{\sqrt{M}|d|}{2\pi} \right) - \frac{\Gamma'}{\Gamma}(1
+ \alpha) - \frac{\Gamma'}{\Gamma}(1 - \alpha).
\end{equation}
For the integral on the $(1-c)$ line we change variables $s
\rightarrow 1-s$ and use (\ref{zwei}). We thus obtain finally the
following:
\begin{theorem} Assuming the Ratios Conjecture \ref{conj:ratios},
the 1-level density for the zeros of the family of even
quadratic twists of an elliptic curve $L$-function $L_E(s)$ with
prime conductor $M$ is given by
\begin{align} \label{oneleveldensity}\nonumber
S_1(f) = &~\sum_{\substack{0<d \leq X\\
\chi_d(-M)\omega_E=+1}} \sum_{\gamma_d} f(\gamma_d)\\=&~
\frac{1}{2\pi} \int_{-\infty}^\infty f(t) \sum_{\substack{0<d \leq X\\
\chi_d(-M)\omega_E=+1}} \Bigg( 2\log
\left(\frac{\sqrt{M}|d|}{2\pi} \right) + \frac{\Gamma'}{\Gamma}(1
+ it) + \frac{\Gamma'}{\Gamma}(1 - it)\nonumber
\\  & + 2\Big[-\frac{\zeta'(1 + 2it)}{\zeta(1 + 2it)} +
\frac{L_E'({\rm sym}^2, 1 + 2it)}{L_E({\rm sym}^2, 1+2it)} +
A_E^1(it,it)
\\ \nonumber & - \left(\frac{\sqrt{M}|d|}{2\pi}\right)^{-2it} \frac{\Gamma(1 - it)}{\Gamma(1 + it)}
\frac{\zeta(1 + 2it)L_E({\rm sym}^2, 1 -
2it)}{L_E({\rm sym}^2,1)}A_E(-it, it)\Big] \Bigg)dt \\
 & + O(X^{1/2+\varepsilon}),\nonumber
\end{align}
where $\gamma_d$ is a generic zero of $L_E(s,\chi_d)$, $f$ is an
even test function as described above, $\omega_E$ is the sign from
the functional equation of $L_E$, $L_E({\rm sym}^2,s)$ is the
associated symmetric square $L$-function (defined at
(\ref{Iwaniec})), and $A_E$ and $A_E^1$ are arithmetic factors
defined at (\ref{AE}) and (\ref{eq:AE1}), respectively.
\end{theorem}


\section{Numerical test} \label{sect:data}
\begin{figure}[p]
\includegraphics[scale=.78]{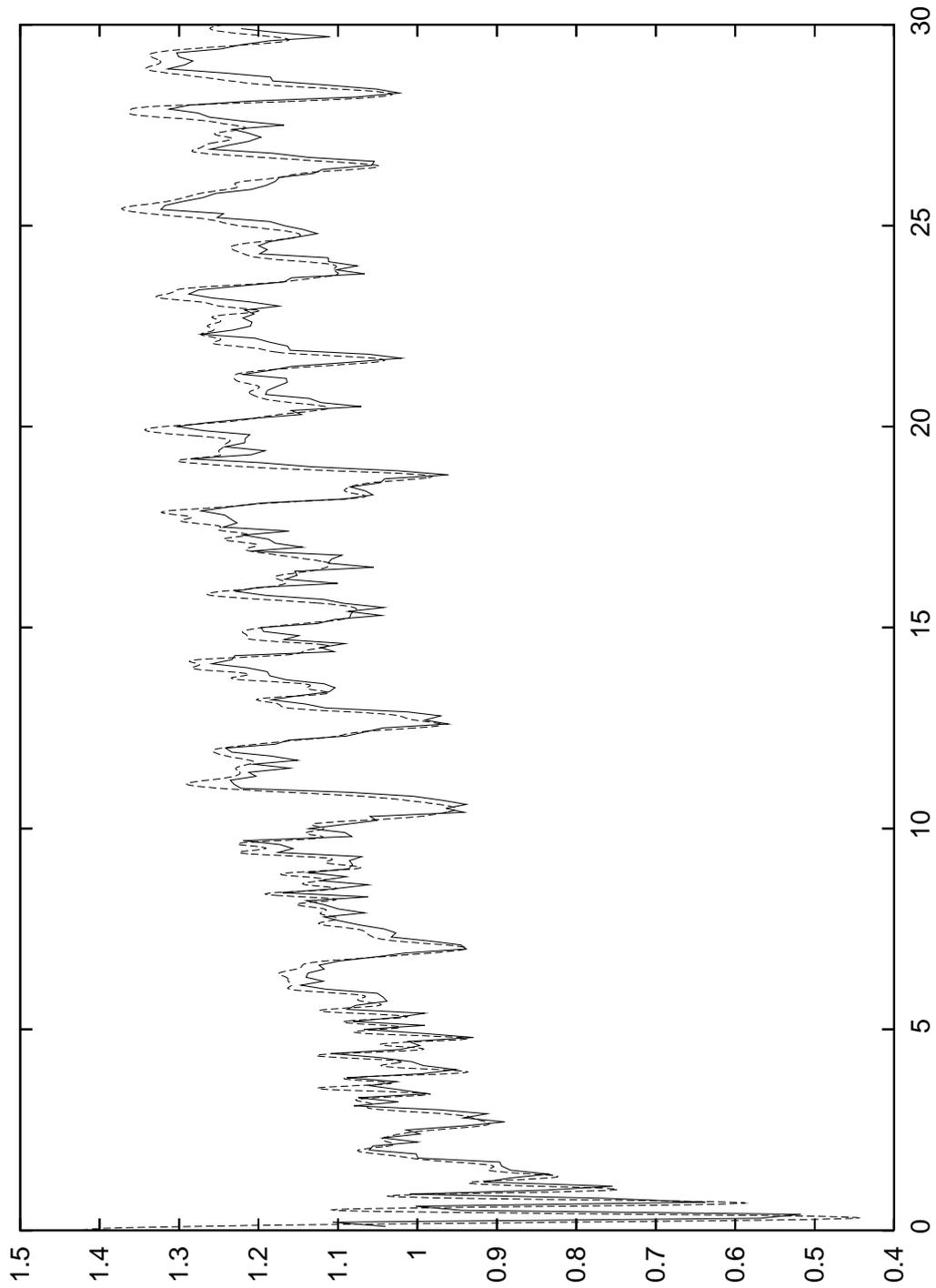}
\caption{1-level density of unscaled zeros from 0 up to height 30
of even quadratic twists of $L_{E_{11}}$ with $0 < d < 40,000$:
prediction (dashed), from (\ref{oneleveldensity}), versus
numerical data (solid)} \label{fig:datavstheory}\end{figure}

We test our prediction -- namely formula (\ref{oneleveldensity})
-- for the 1-level density with a concrete example (see figure
\ref{fig:datavstheory}).  We pick the elliptic curve $E_{11}$ with
$(a_1,a_2,a_3,a_4,a_6) = (0,-1,1,0,0)$ in the Weierstraß form
\begin{equation}
y^2 + a_1 xy + a_3 y = x^3 + a_2 x^2 + a_4 x + a_6
\end{equation}
giving
\begin{equation}
E_{11}: y^2 + y = x^3 - x^2
\end{equation}
and consider the even quadratic twists of its associated
$L$-function with fundamental discriminants between 0 and 40,000.
We are interested in the 1-level density of unscaled zeros from 0
up to height 30. The numerical data is obtained from Rubinstein's
{\sffamily lcalc} \cite{kn:rubinstein}. In the range considered we
find 11,135 quadratic twists, of which 5,562 are even ones with a
total of about 590,170 zeros.  In figure \ref{fig:datavstheory} we
obtain the solid curve from the histogram of this zero data by
choosing a binsize of 0.1 and dividing by both the number of
quadratic twists with even functional equation, and the mean
density of zeros $\log(\sqrt{11}X/(2\pi))$.  593 of the
$L$-functions with even functional equation have (at least) a
double zero at the central point; these zeros at the central point
are not plotted in figure \ref{fig:datavstheory}.  The dashed
curve is obtained from the formula (\ref{oneleveldensity}) with
$X=40,000$ and $f(t)=\delta(t-x)+\delta(t+x)$ for $x$ between 0
and 30.  This curve is scaled like the data curve by dividing
through by the number of quadratic twists with even functional
equation and the mean density of zeros.  It was computed using a
combination of Mathematica and C++. The coefficients $\lambda(p)$
appearing in the arithmetic factor $A_E(\alpha, \gamma)$ were
computed using PARI.  To compute coefficients of prime powers
$\lambda(p^m)$ for $p \nmid M$ the following recursion formulas
(see \cite{kn:hugmil07}) were used
\begin{eqnarray}
\lambda(p^{2m}) & = & \lambda(p)^{2m} - \sum_{r = 0}^{m - 1} \left({2m \choose m - r} - {2m \choose m - r - 1} \right)\lambda(p^{2r})\\
\lambda(p^{2m + 1}) & = & \lambda(p)^{2m + 1} - \sum_{r = 0}^{m - 1} \left({2m + 1 \choose m - r}
- {2m + 1 \choose m - r - 1} \right)\lambda(p^{2r + 1}).
\end{eqnarray}

In general there is good agreement between the data and the
theoretical curve, which captures the main features of the data.
We would expect better agreement with a larger set of data, since
the data seems not yet to have resolved all the peaks further out
along the axis.

\begin{figure}[p]
\includegraphics[scale=.7]{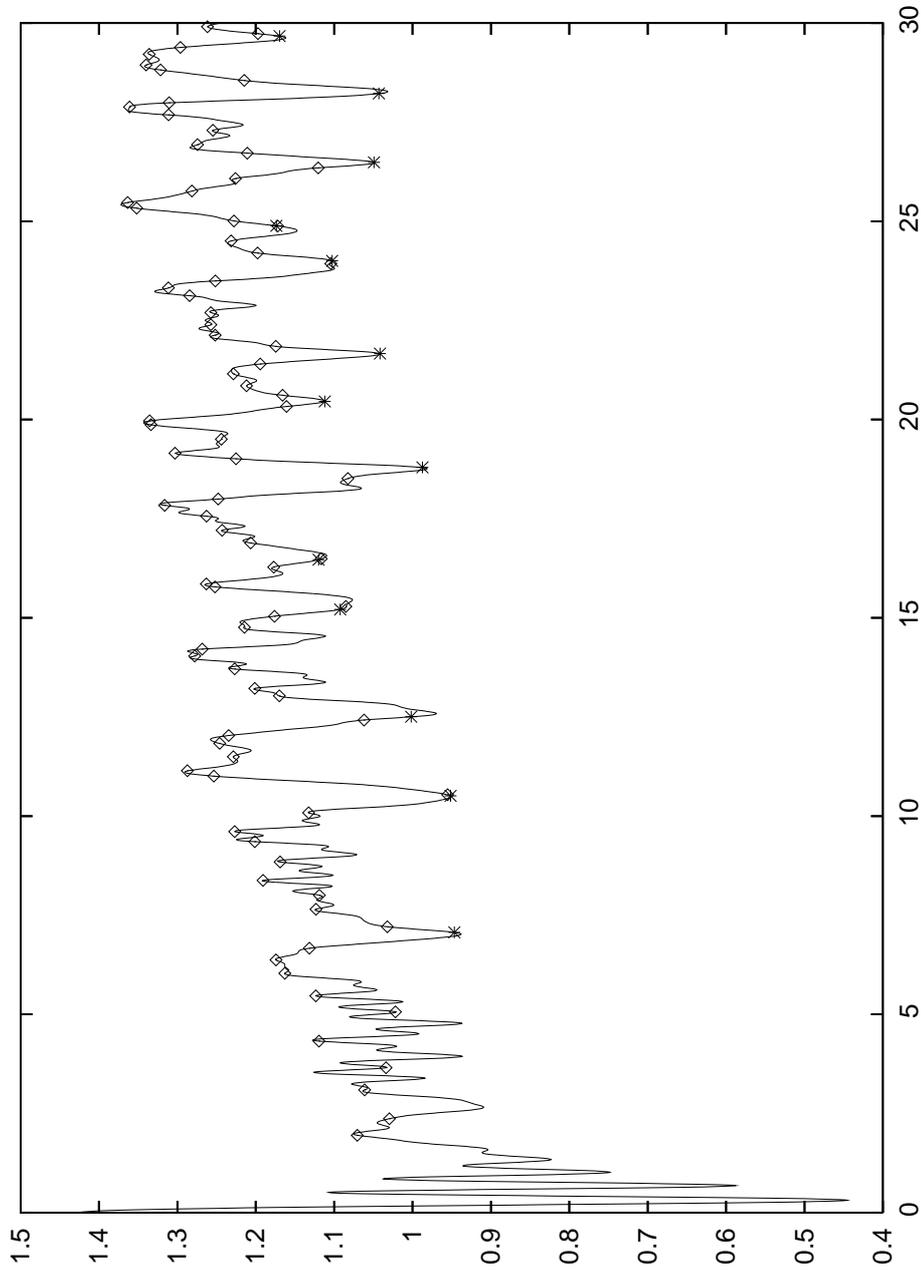}
\caption{Effects of non-trivial zeros of the Riemann zeta function
(indicated by~$\ast$)  and the non-trivial zeros of $L(\sym^2, s)$
function (indicated by $\diamond$) on the conjectural formula
(\ref{oneleveldensity}) for the 1-level density of unscaled zeros
from 0 up to height 30 of even quadratic twists of $L_{E_{11}}$
with $0 < d < 40,000$}\label{fig:zetaLzeros}
\end{figure}

A closer look reveals that the 1-level-density is strongly
governed by the non-trivial zeros of $\zeta(s)$ and $L(\sym^2,
s)$: we observe that some dips of the data curve are located at
$\gamma / 2$ where $\gamma$ is the ordinate of a non-trivial zero
of the Riemann zeta function. This is captured in the term
\begin{equation}
-\frac{\zeta'(1+2it)}{\zeta(1+2it)}
\end{equation}
of our conjecture for $S_1(f)$. In figure \ref{fig:zetaLzeros} we
mark the position of a non-trivial zero of the Riemann zeta
function on our conjectural answer by a $\ast$. These $\ast$ are
all localised in or around a neighbourhood of a dip.  This
phenomenon has been encountered before, in the study of lower order terms
of the number variance \cite{kn:berry88} and the correlation
functions
\cite{kn:berrykeating99,kn:bk96,kn:consna06,kn:consna07,kn:consna08}
of the Riemann zeros, and in the one-level density of other
families of $L$-functions \cite{kn:consna06}.


On the other hand we observe that some peaks are located at
$\tilde{\gamma} / 2$ where $\tilde{\gamma}$ is the ordinate of a
non-trivial zero of $L_E(\sym^2, s)$. This is captured in the term
\begin{equation}
\frac{L_E'(\sym^2, 1 + 2it)}{L_E(\sym^2, 1 + 2it)}
\end{equation}
of our conjecture for $S_1(f)$. In figure \ref{fig:zetaLzeros} we
mark the position of a non-trivial zero of $L_E(\sym^2, s)$ by a
$\diamond$. The majority of these $\diamond$s are localized in or
around a neighbourhood of a peak.  In particular, we observe that
if a zero of the Riemann zeta function is close to a zero of
$L(\sym^2, s)$ then these zeros are localised in or around a dip.
Hence, zeros of the Riemann zeta function seem to dominate the
behaviour of the 1-level-density more than the zeros of $L(\sym^2,
s)$. This may be explained because the density of the Riemann
zeros in this range is smaller than that of the zeros of
$L(\sym^2, s)$ and so in terms of the mean zero density the
one-line is closer to the half-line in the case of the Riemann
zeta function.  Therefore one would expect the Riemann zeros to
have a larger effect.

The term
\begin{equation}
- \left(\frac{\sqrt{M}|d|}{2\pi}\right)^{-2it} \frac{\Gamma(1 -
it)}{\Gamma(1 + it)} \frac{\zeta(1 + 2it)L_E(\sym^2, 1 -
2it)}{L_E(\sym^2,1)}A_E(-it, it),
\end{equation}
from (\ref{oneleveldensity}), makes its most obvious contribution
by causing the oscillation near the origin of the plot of our
conjectural answer for the 1-level density.  The factor
$\left(\frac{\sqrt{M}|d|}{2\pi}\right)^{-2it}$ results in
oscillations on the scale of the mean density of the zeros of the
original $L$-function, $L_E$.

\begin{figure}[htbp]
\begin{center}
\includegraphics[scale=.27,
angle=-90]{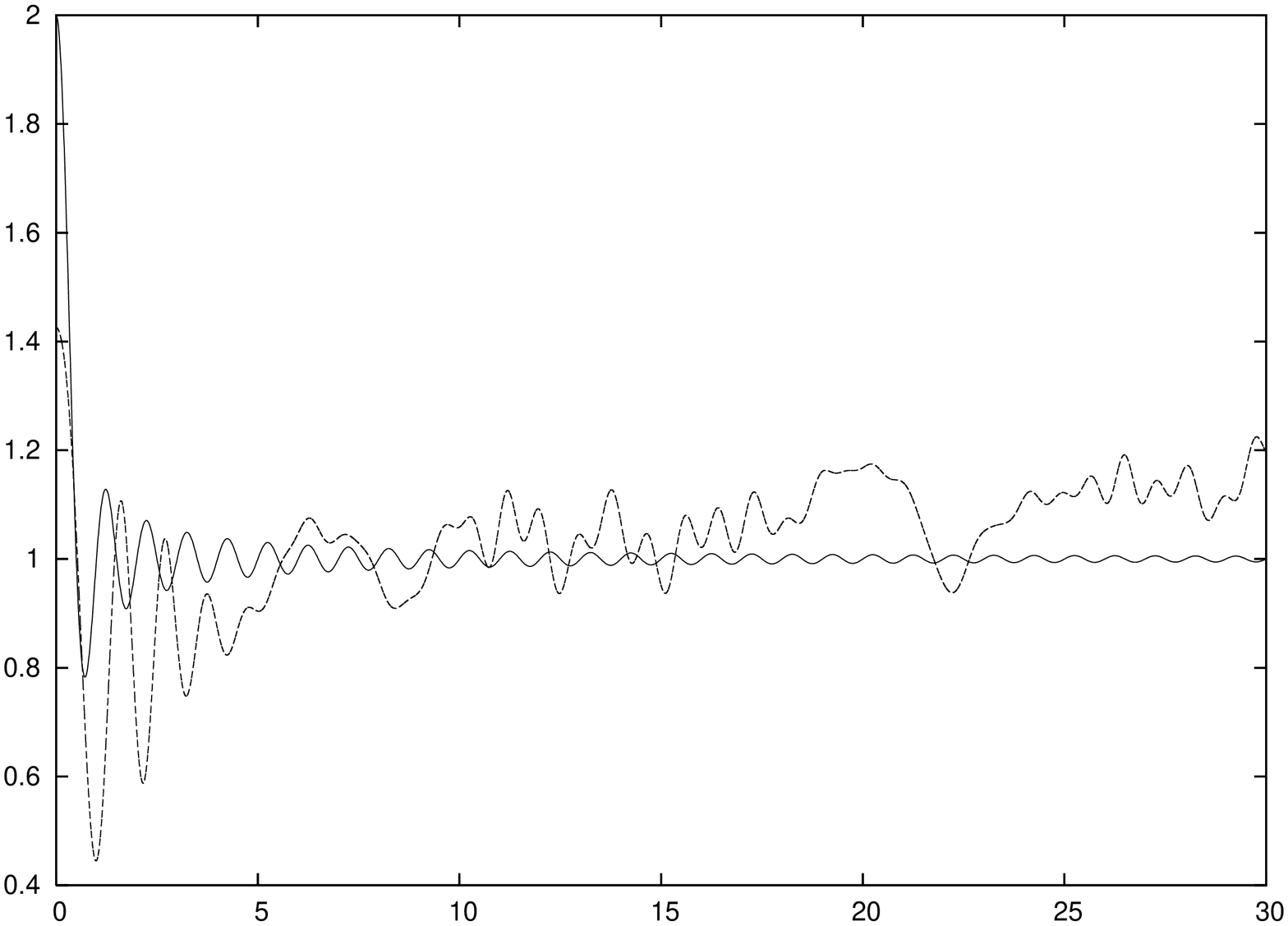}\hspace{0.05
in}
\includegraphics[scale=.27,
angle=-90]{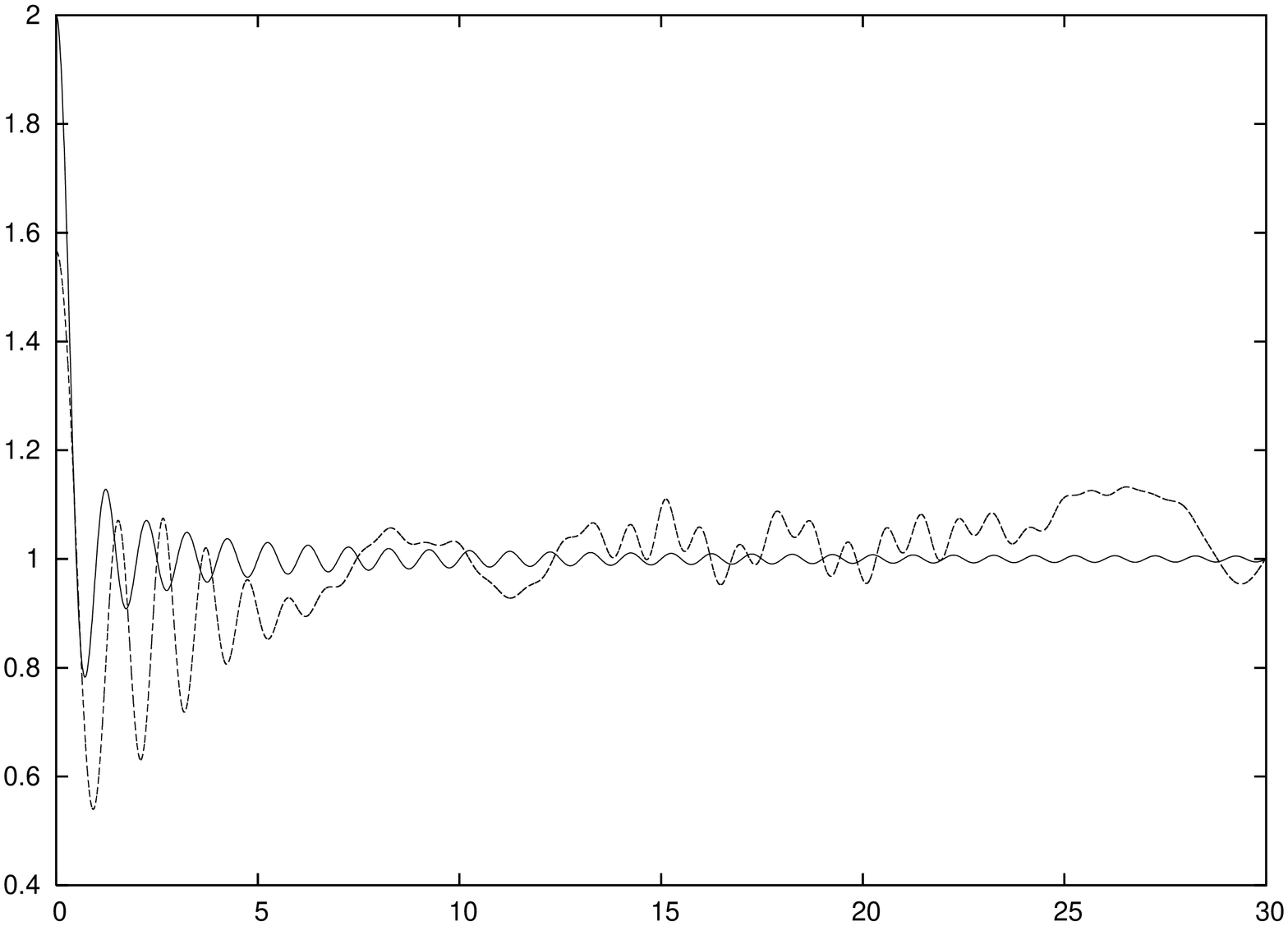}\\
\includegraphics[scale=.27,
angle=-90]{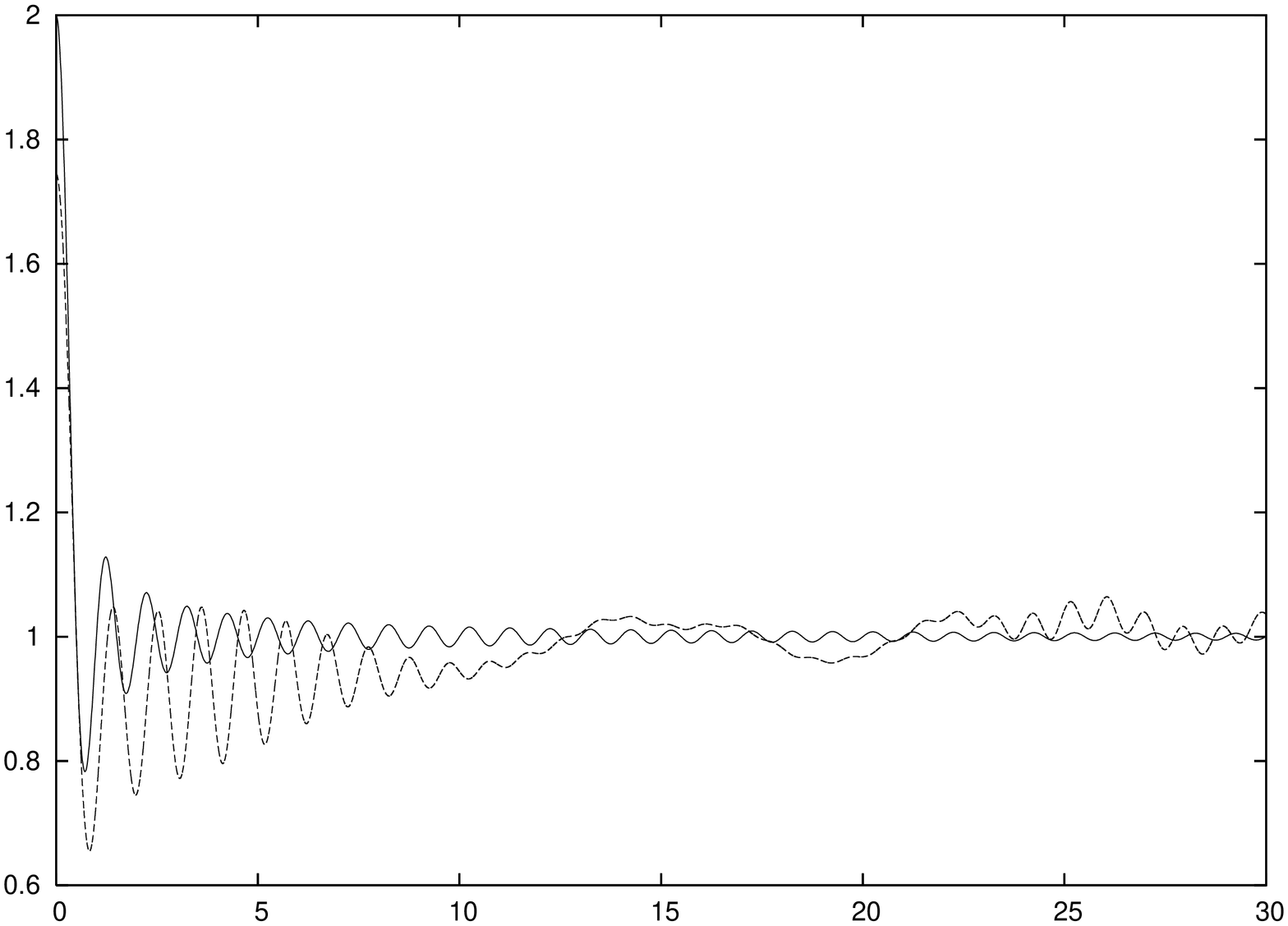}\hspace{0.05
in}\includegraphics[scale=.27,angle=-90]{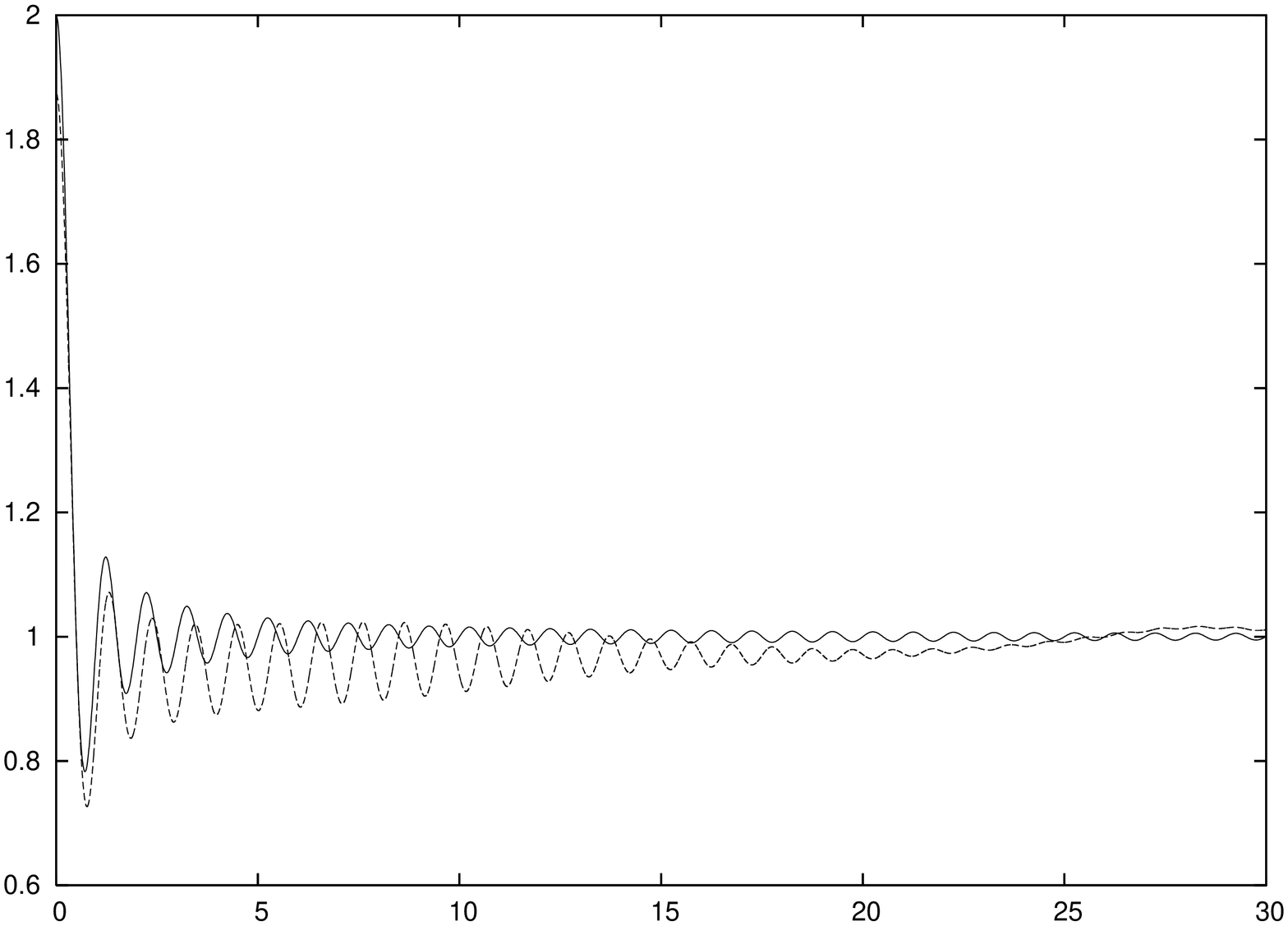}\\
\includegraphics[scale=.27,
angle=-90]{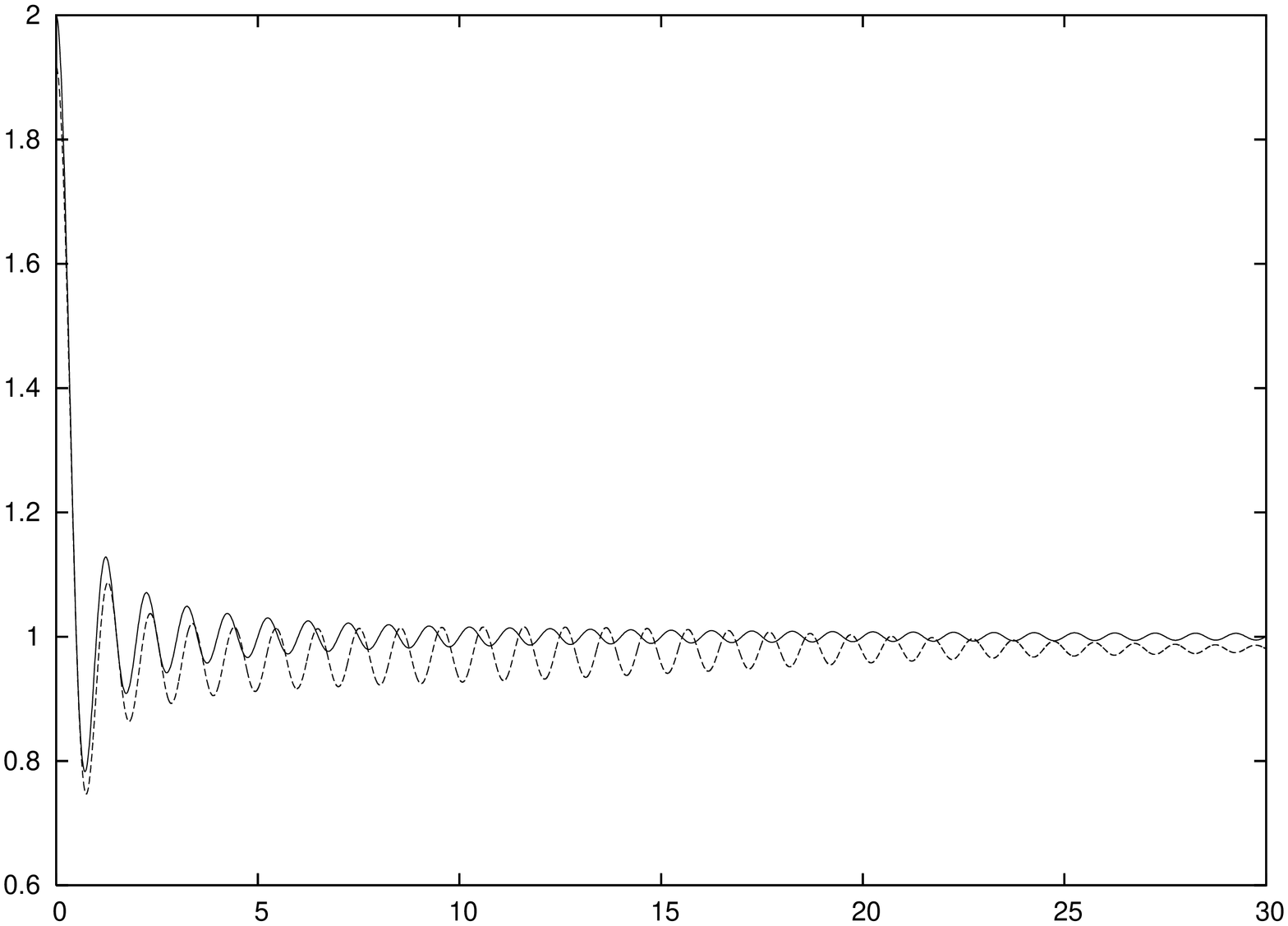}\hspace{0.05
in}
\includegraphics[scale=.27,
angle=-90]{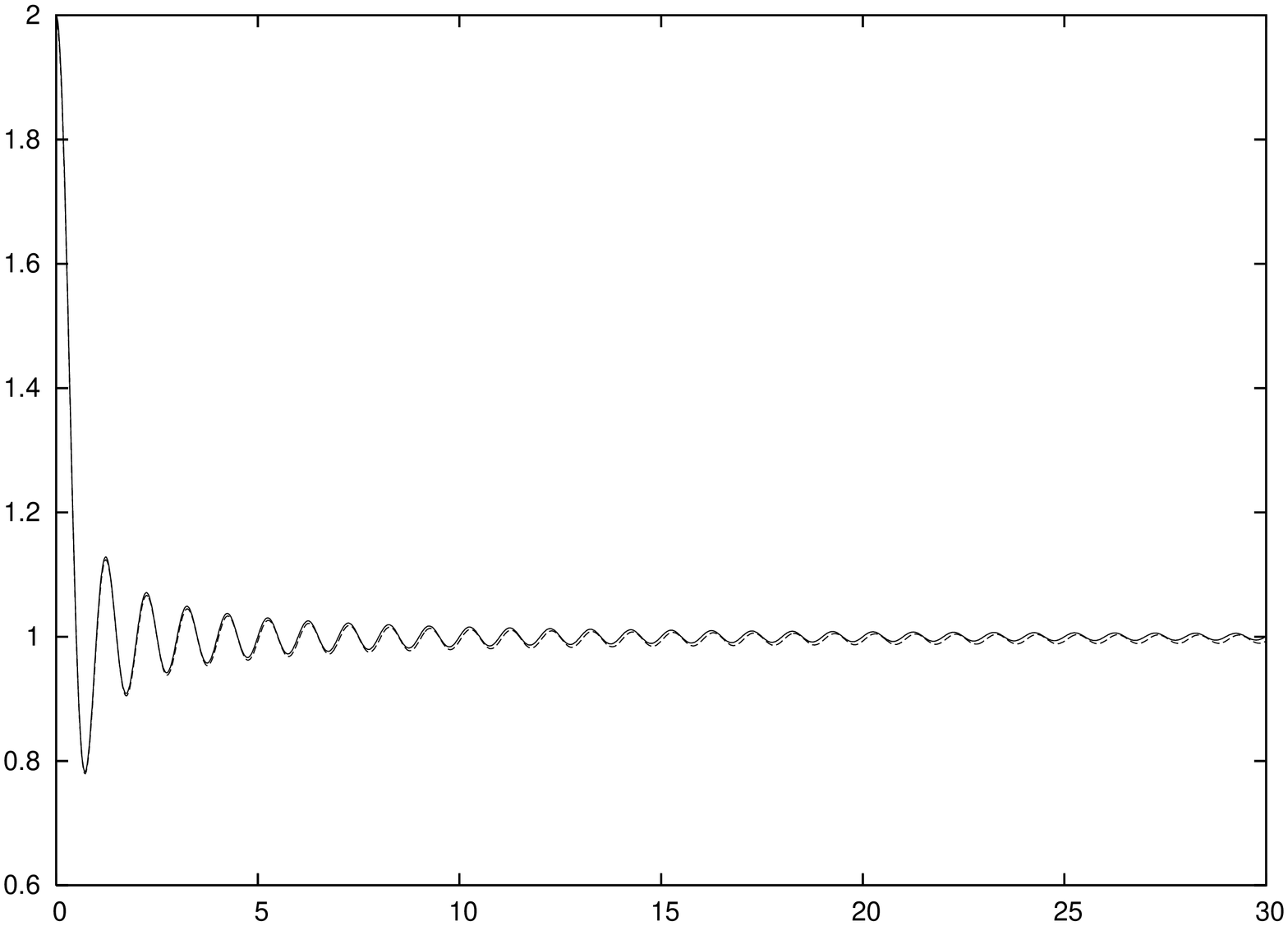}
\caption{Scaled limiting 1-level density of $SO(2N)$ (solid)
versus scaled  formula (\ref{interim}) divided by $X^*$ (dashed)
for: $0 < d \leq 40,000$ (top left), $0 < d \leq 10^6$ (top
right), $0 < d \leq 10^{10}$ (middle left), $0 < d \leq 10^{20}$
(middle right), $0 < d \leq 10^{30}$ (bottom left), $0 < d \leq
10^{300}$ (bottom right)} \label{fig:scale1}
\end{center}
\end{figure}

In summary, we notice that the lower order terms dominate the
behaviour of the zeros when we are far from the limit of infinite
conductor (in the family of quadratic twists, $E_d$, the conductor
increases with $d$). This becomes more obvious when we compare our
conjectural answer for finite conductors with the limiting
theoretical result: in figure \ref{fig:scale1}  we consider the
scaled 1-level density of $SO(2N)$ in the limit $N\rightarrow
\infty$ against our conjectural answer (also scaled) for finite
conductor. We observe convergence to the limiting theoretical
result as we increase $X$, the cut-off point for $d$. The observed
effects of the arithmetical terms for small and finite conductors
are washed out and shifted away from the origin in the large
conductor limit.

To further understand the approach to the limiting distribution,
we calculate the 1-level density for scaled zeros and recover
the limit and the next to leading order term from
(\ref{oneleveldensity}). As a first step we rescale the variable
$t$ in (\ref{oneleveldensity}) as
\begin{equation}
\tau = t (L/\pi)
\end{equation}
and define
\begin{equation}\label{eq:fg}
f(t) = g(t (L/\pi)),
\end{equation}
where
\begin{equation}
L := \log\bigg(\frac{\sqrt{M}X}{2\pi}\bigg),
\end{equation}
and get, after a change of variables,
\begin{eqnarray}
&& \sum_{\substack{0<d \leq X\\
\chi_d(-M)\omega_E=+1}} \sum_{\gamma_d}g\Big(\frac{\gamma_d
L}{\pi}\Big)\nonumber
\\ \nonumber && =~ \frac{1}{2L} \int_{-\infty}^\infty g(\tau) \sum_{\substack{0<d \leq X\\
\chi_d(-M)\omega_E=+1}} \Bigg( 2\log
\left(\frac{\sqrt{M}|d|}{2\pi} \right) +
\frac{\Gamma'}{\Gamma}\Big(1 + \frac{i\pi \tau}{L}\Big)
\\ \nonumber &&~~~ + \frac{\Gamma'}{\Gamma}\Big(1 - \frac{i \pi \tau}{L}\Big) +
2\Big[-\frac{\zeta'(1 + \frac{2 i \pi \tau}{L})}{\zeta(1 + \frac{2
i \pi \tau}{L})} + \frac{L_E'(\sym^2, 1 + \frac{2 i \pi
\tau}{L})}{L_E(\sym^2, 1+\frac{2 i \pi \tau}{L})} +
A_E^1\Big(\frac{i \pi \tau}{L}, \frac{i \pi \tau}{L}\Big)
\\ \nonumber &&~~~ -\bigg(\frac{\sqrt{M}|d|}{2\pi}\bigg)^{-2i \pi \tau / L}
\frac{\Gamma(1 - \frac{i \pi \tau}{L})}{\Gamma(1 + \frac{i \pi
\tau}{L})} \frac{\zeta(1 + \frac{2 i \pi \tau}{L})L_E(\sym^2, 1 -
\frac{2 i \pi \tau}{L})}{L_E(\sym^2,1)}\\
\nonumber&&~~~\qquad\qquad\qquad\qquad\qquad\qquad\qquad\qquad\qquad\qquad
\times A_E\Big(-\frac{i \pi
\tau}{L}, \frac{i \pi \tau}{L}\Big)\Big] \Bigg)d\tau \\
\label{interim} &&~~~ + O(X^{1/2+\varepsilon}).
\end{eqnarray}

We write the number of fundamental discriminants less than or
equal to $X$ as
\begin{equation}
X^* := \sum_{\substack{0<d \leq X\\
\chi_d(-M)\omega_E=+1}} 1.
\end{equation}
Using the Euler-Maclaurin formula we make the approximation
\begin{equation}
\sum_{\substack{0<d \leq X\\
\chi_d(-M)\omega_E=+1}} \log \left(\frac{\sqrt{M}|d|}{2\pi}
\right) = X^*\left[ \log \left(\frac{\sqrt{M}X}{2\pi} \right) - 1
\right] + O\left(X^{1/2 + \varepsilon}\right).
\end{equation}
In the same manner we have
\begin{align}
\sum_{\substack{0<d \leq X\\
\chi_d(-M)\omega_E=+1}} \left(\frac{\sqrt{M}|d|}{2
\pi}\right)^{-2i \pi \tau / L} = & X^*\big(1+ \frac{2i \pi
\tau}{L} +O(L^{-2})\big)e^{-2i\pi\tau} + O (X^{1/2}).
\kommentar{
\frac{X^*}{1- (2i \pi \tau)/L }
\bigg(\frac{\sqrt{M}X}{2\pi} \bigg)^{-2i \pi \tau/L} +
O(X^{-1/2})\nonumber \\
& =& ******X^*\big(1+ \frac{2i \pi \tau}{L} +O(\log^{-2}
L)\big)e^{-2i\pi\tau}***** .}
\end{align}
Writing
\begin{equation}
 \zeta(s+1) = \frac{1}{s} + \sum_{n=0}^\infty
\frac{(-1)^n}{n!} \gamma_n s^n,
\end{equation}
we have
\begin{align}
\frac{\zeta'(1+s)}{\zeta(1+s)} = & -s^{-1} + \gamma + (-\gamma^2 -
2\gamma_1)s  + O(s^2),
\end{align}
where $\gamma = \gamma_0$ is Euler's constant, and so
\begin{equation}
\zeta(1 + \frac{2 i \pi \tau}{L}) = \frac{L}{2 i \pi \tau} + \gamma + O(L^{-1}).
\end{equation}
 and
\begin{equation}
\frac{\zeta'(1 + \frac{2 i \pi \tau}{L})}{\zeta(1 + \frac{2 i \pi
\tau}{L})}=-\frac{L}{2i \pi \tau} + \gamma +O(L^{-1}).
\end{equation}
Simple Taylor expansions of the other factors in (\ref{interim})
lead us to, with the relation between $f$ and $g$ given in
(\ref{eq:fg}),
\begin{eqnarray} \label{S1f}&& \frac{1}{X^*}S_1(f)=\frac{1}{X^*}\sum_{\substack{0<d \leq X\\
\chi_d(-M)\omega_E=+1}} \sum_{\gamma_d}g\Big(\frac{\gamma_d L}{\pi}\Big)\\
&&\quad=\int_{-\infty}^{\infty}g(\tau)\Bigg( 1 + \frac{\sin(2\pi
\tau)}{2 \pi \tau} - a_1\frac{1 +\cos(2\pi \tau)}{L} -
a_2\frac{\pi \tau \sin(2\pi \tau)}{L^2} +
O\left(\frac{1}{L^3}\right)\Bigg)d\tau\nonumber
\end{eqnarray}
where
\begin{equation} \label{aone}
a_1 = 1 + 2\gamma - A^1_E(0,0) - \frac{L'_E(\sym^2,
1)}{L_E(\sym^2, 1)}
\end{equation}
and
\begin{align} \nonumber
a_2 & = 2 + 4\gamma + 3\gamma^2 - 2\gamma_1 + B'(0) + 2\gamma B'(0) - 2\frac{L'_E(\sym^2, 1)}{L_E(\sym^2, 1)}\\
& -\frac{4\gamma L'(1)}{L(1)} - \frac{B'(0)L'_E(\sym^2,
1)}{L_E(\sym^2, 1)} + \frac{B''(0)}{4} + \frac{L''_E(\sym^2,
1)}{L_E(\sym^2, 1)},
\end{align}
with
\begin{equation}
B'(0) = \frac{d}{dr}A_E(-r,r)\Big|_{r =0} \mbox{~and~} B''(0) =
\frac{d^2}{dr^2}A_E(-r,r)\Big|_{r =0}.
\end{equation}
In order to obtain (\ref{aone}) we use the following identity
\begin{equation} \label{relation}
-\frac{1}{2} B'(0) = A^1_E(0,0).
\end{equation}
We establish identity (\ref{relation}) by simple algebra and using
(\ref{lang}) for primes not dividing $M$, the multiplicativity of
$\lambda(p)$ for $p|M$ and $A_E(r,r) =1$.

This work was initially conceived to investigate the unexpected
numerical results found by Steven J. Miller \cite{kn:mil05} near
the origin of the histogram of the distribution of the first zero
above the central point of a family of rank zero $L$-functions. He
observed very few examples of zeros lying close to the central
point.  That is, he observed the phenomenon of repulsion of zeros
from the central point which we know, from rigorous work on the
1-level and 2-level densities \cite{kn:young,kn:miller02,kn:mil04}
does not persist in the large conductor limit.  Since the 1-level
density (a histogram of all zeros) and the distribution of the
lowest zero (a histogram of the lowest zero of each $L$-function)
are the same for very small distances from the central point, it
is natural to enquire whether the ratios conjecture yields a
formula for the 1-level density which would display and explain
Miller's observed repulsion at finite conductor. Although it can
be seen from figure \ref{fig:scale1} that the formula
(\ref{interim}) is significantly smaller near the origin than the
limiting curve, and approaches it from below as the conductor
increases, there is no evidence of repulsion. This is a major
discrepancy from the data, as seen in figure~\ref{fig:nrzero}:
away from the critical point we have a nice match between the
prediction and the data while near the critical point we find
fewer zeros in the data than predicted by our formula. It is most
interesting that the main terms of the ratios conjecture do not
capture this important feature.  Of course, the natural question
is whether this contradicts the ratios conjecture, or whether the
discrepancy can by accounted for by the error term. As expected
due to the limited data available, the test described below is
inconclusive, but shows signs that the error term in the
ratios conjecture (and hence on the one level density in
(\ref{oneleveldensity})) is of the form $X^{b+\varepsilon}$, for
$b<1$.  The ratios conjecture is usually stated with $b=1/2$.

We fix several sample points at various distances away from the
critical point and measure the difference between the main terms
of our prediction  (that is, the sum over $d$ inside the integral
in (\ref{oneleveldensity})) and the data. In fact, we compare the
normalised versions of our prediction and data by dividing through
by the number of fundamental discriminants $X^*$ less than $X$ and
the mean density of zeros. So let us denote this difference
between the main terms of the normalised theory and the data at a
fixed height $t$ and fixed $X$ by $\Delta(t, X)$.  Since we have
divided by $X^*$, which is proportional to $X$, this difference is
expected to be of size
\begin{equation}
|\Delta(t, X)| = O(X^{b-1 + \varepsilon})
\end{equation}

The quantity we will plot is
\begin{equation} \label{Q_difference}
Q_{\Delta}(t,X) := \frac{\log(|\Delta( t, X)|)}{\log X}
\end{equation}
and if the ratios conjecture with error term $X^{b+\varepsilon}$
is correct then we would expect
\begin{equation} \label{log_difference}
Q_{\Delta}(t,X) = b-1 + O\Big(\frac{\log\log X}{\log X}\Big)
\end{equation}
as $X\rightarrow \infty$.

\begin{figure}[h]
\includegraphics[scale=.55,angle=-90]{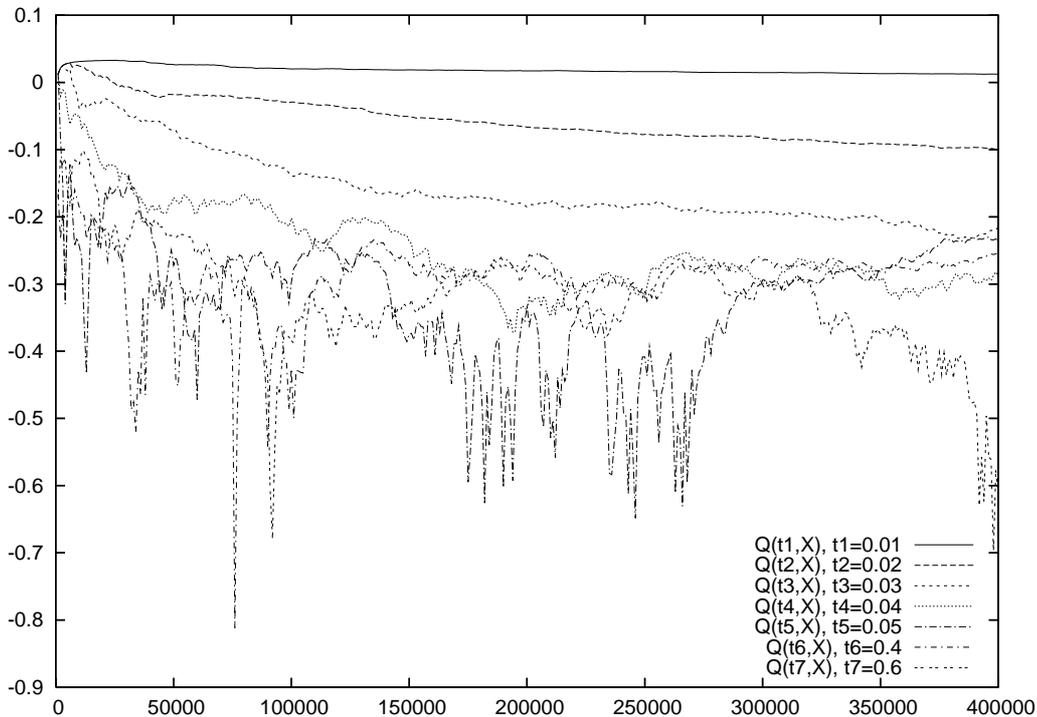}
\caption{discrepancy $Q_{\Delta}(t,X)$ from (\ref{Q_difference}) between prediction and
data for the 1-level density of even quadratic twists of $L_{E_{11}}$ with $0 < d < 400,000$.
}
\label{fig:discrepancy}
\end{figure}

In figure \ref{fig:discrepancy} we plot the quantity
$Q_{\Delta}(t,X)$ for $0 < X < 400,000$ and for various fixed
sample points $ t_1~=~0.01, t_2~=~0.02, t_3~=~0.03, t_4~=~0.04,
t_5~=~0.05, t_6~=~0.4$ and, $t_7~=~0.6$.  We notice that the
curves are much smoother for sample points near the critical
point, $t=0$, eg.~ $t_1,t_2,t_3$.  In the range $0 < X < 400,000$
these points are well inside the region where the zero data shows
repulsion at the critical point; see figure \ref{fig:nrzero}. Thus
the difference between the theory (smooth curve in figure
\ref{fig:nrzero}) and data (histogram) does not change sign as $X$
increases.  Presumably it is the amplification of such sign
changes by the logarithm in (\ref{Q_difference}) that is
responsible for the jagged curves in figure \ref{fig:discrepancy}
for sample points $t_4, t_5$ and $t_6$.

We see also that the curves at sample points close to the critical
point appear at first sight to indicate a larger error term - in
fact, over this range of $X$ the $t_1$ curve implies $b-1>0$!  If
a limit such as (\ref{log_difference}) exists, it does not seem to
behave uniformly in $t$.  However, the $t_1, t_2$ and $t_3$ curves
are decaying as $X$ increases and we do not have enough data to
see what their final behaviour will be.  We remember that the
convergence is like $\log \log X/ \log X$, so we would need much
more data to be able to make a sensible conclusion about the size
of the error term.

Also, it is interesting to note that at the right hand side of
figure \ref{fig:discrepancy} the $t_3=0.03$ curve has decayed to a
level comparable to the curves of the sample points that are more
distant from $t=0$.  Examining figure \ref{fig:nrzero}, it appears
that the area of major discrepancy between the ratios conjecture
prediction and the data (that is, where the data shows repulsion
from the critical point at $t=0$) lies between t=0 and about
$t=0.03$.  We expect that this region will narrow as the range of
discriminants, $d$, increases, and this is born out by comparing
the two pictures in figure \ref{fig:nrzero}; the data grows more
quickly to the height of the solid curve in the right hand picture
where $0<d<400,000$, than in the left hand picture where
$0<d<100,000$. Thus at the right hand edge of figure
\ref{fig:discrepancy}, the point $t_3=0.03$ is about to move into
the region where there is good agreement between the ratios
conjecture prediction and the data. Making a speculative
conclusion from the limited data available, this suggests that the
curves for $t_1$ and $t_2$, or any other fixed $t$, would also
decay to this level if we could gather enough data to shrink the
area of discrepancy at the origin of figure~\ref{fig:nrzero} to a
narrow enough band.

It is impossible to say from the available data what the exponent
$b$ in the error term of the ratios conjecture is.  There is
certainly  no evidence to suggest $b=0.5$, but the possibility
that the curves in figure \ref{fig:discrepancy} would decay to
$-0.5$ if we could vastly extend the rage of the plot is not ruled
out. However, figure \ref{fig:discrepancy} certainly appears to
suggest that $b<0$ and so the error term is a power of $X$ smaller
than the main term.

\section{ Summary}

We find that the ratios  conjecture provides a formula for the one
level density of zeros of a family of quadratic twists of an
elliptic curve $L$-function that agrees with data for finite
conductor, except in the vicinity of the critical point, $t=0$,
and explains the arithmetic nature of the lower order terms which
entirely dominate the behaviour of the statistic away from $t=0$.
The ratios conjecture prediction, when properly scaled, approaches
the limiting $SO(2N)$ random matrix result as the family of
elliptic curves includes those with larger and larger conductor.
This supports all the available evidence that $SO(2N)$ is the
correct limit for zero statistics in this family. It is very
interesting that the ratios conjecture prediction does not capture
the phenomenon of zero repulsion from the critical point, $t=0$,
but the data we have available certainly allows for the ratios
conjecture to be correct with some power $b<1$ of $X$ in the error
term; the discrepancy between the ratios conjecture prediction and
the data (at the origin of figure \ref{fig:nrzero}) can quite
possibly be contained in the error term.

In ongoing work of the authors in collaboration with E. Due{\~n}ez
and S. J. Miller we propose an explanation for the observed
repulsion of zeros near the central point for finite conductor and
a random matrix model that captures the phenomenon.


\pagebreak

\newcommand{\etalchar}[1]{$^{#1}$}


\end{document}